\documentclass[11pt]{article}

\usepackage{amssymb,latexsym,amsmath}

\usepackage{graphicx}

\begin{document}

\newcommand{\bfi}{\bfseries\itshape}

\makeatletter

\@addtoreset{figure}{section}

\def\thefigure{\thesection.\@arabic\c@figure}

\def\fps@figure{h, t}

\@addtoreset{table}{bsection}

\def\thetable{\thesection.\@arabic\c@table}

\def\fps@table{h, t}

\@addtoreset{equation}{section}

\def\theequation{\thesubsection.\arabic{equation}}

\makeatother

\newtheorem{thm}{Theorem}[section]

\newtheorem{prop}[thm]{Proposition}

\newtheorem{lema}[thm]{Lemma}

\newtheorem{cor}[thm]{Corollary}

\newtheorem{defi}[thm]{Definition}

\newtheorem{rk}[thm]{Remark}

\newtheorem{exempl}{Example}[section]

\newenvironment{exemplu}{\begin{exempl}  \em}{\hfill $\square$

\end{exempl}}

\newcommand{\comment}[1]{\par\noindent{\raggedright\texttt{#1}

\par\marginpar{\textsc{Comment}}}}

\newcommand{\todo}[1]{\vspace{5 mm}\par \noindent \marginpar{\textsc{ToDo}}\framebox{\begin{minipage}[c]{0.95 \textwidth}

\tt #1 \end{minipage}}\vspace{5 mm}\par}

\newcommand{\ea}{\mbox{{\bf a}}}

\newcommand{\eu}{\mbox{{\bf u}}}

\newcommand{\ueu}{\underline{\eu}}

\newcommand{\ueo}{\overline{u}}

\newcommand{\oeu}{\overline{\eu}}

\newcommand{\ew}{\mbox{{\bf w}}}

\newcommand{\ef}{\mbox{{\bf f}}}

\newcommand{\eF}{\mbox{{\bf F}}}

\newcommand{\eC}{\mbox{{\bf C}}}

\newcommand{\en}{\mbox{{\bf n}}}

\newcommand{\eT}{\mbox{{\bf T}}}

\newcommand{\eL}{\mbox{{\bf L}}}

\newcommand{\eR}{\mbox{{\bf R}}}

\newcommand{\eV}{\mbox{{\bf V}}}

\newcommand{\eU}{\mbox{{\bf U}}}

\newcommand{\ev}{\mbox{{\bf v}}}

\newcommand{\eve}{\mbox{{\bf e}}}

\newcommand{\uev}{\underline{\ev}}

\newcommand{\eY}{\mbox{{\bf Y}}}

\newcommand{\eK}{\mbox{{\bf K}}}

\newcommand{\eP}{\mbox{{\bf P}}}

\newcommand{\eS}{\mbox{{\bf S}}}

\newcommand{\eJ}{\mbox{{\bf J}}}

\newcommand{\eB}{\mbox{{\bf B}}}

\newcommand{\eH}{\mbox{{\bf H}}}

\newcommand{\leb}{\mathcal{ L}^{n}}

\newcommand{\eI}{\mathcal{ I}}

\newcommand{\eE}{\mathcal{ E}}

\newcommand{\hen}{\mathcal{H}^{n-1}}

\newcommand{\eBV}{\mbox{{\bf BV}}}

\newcommand{\eA}{\mbox{{\bf A}}}

\newcommand{\eSBV}{\mbox{{\bf SBV}}}

\newcommand{\eBD}{\mbox{{\bf BD}}}

\newcommand{\eSBD}{\mbox{{\bf SBD}}}

\newcommand{\ecs}{\mbox{{\bf X}}}

\newcommand{\eg}{\mbox{{\bf g}}}

\newcommand{\paromega}{\partial \Omega}

\newcommand{\gau}{\Gamma_{u}}

\newcommand{\gaf}{\Gamma_{f}}

\newcommand{\sig}{{\bf \sigma}}

\newcommand{\gac}{\Gamma_{\mbox{{\bf c}}}}

\newcommand{\deu}{\dot{\eu}}

\newcommand{\dueu}{\underline{\deu}}

\newcommand{\dev}{\dot{\ev}}

\newcommand{\duev}{\underline{\dev}}

\newcommand{\weak}{\stackrel{w}{\approx}}

\newcommand{\mild}{\stackrel{m}{\approx}}

\newcommand{\strong}{\stackrel{s}{\approx}}

\newcommand{\weakdown}{\rightharpoondown}

\newcommand{\opg}{\stackrel{\mathfrak{g}}{\cdot}}

\newcommand{\opunu}{\stackrel{1}{\cdot}}
\newcommand{\opdoi}{\stackrel{2}{\cdot}}

\newcommand{\opn}{\stackrel{\mathfrak{n}}{\cdot}}
\newcommand{\opx}{\stackrel{x}{\cdot}}

\newcommand{\tr}{\ \mbox{tr}}

\newcommand{\Ad}{\ \mbox{Ad}}

\newcommand{\ad}{\ \mbox{ad}}

\renewcommand{\contentsname}{ }

\title{On the 
Kirchheim-Magnani counterexample to metric differentiability}

\author{Marius Buliga \\
\\
Institute of Mathematics, Romanian Academy \\
P.O. BOX 1-764, RO 014700\\
Bucure\c sti, Romania\\
{\footnotesize Marius.Buliga@imar.ro}}

\date{This version:  04.10.2007}

\maketitle


In Kirchheim-Magnani \cite{kirch} the authors construct a left invariant 
 distance $\rho$ on the Heisenberg group such that the identity map $id$ is 
 1-Lipschitz but it is not metrically differentiable anywhere.

In this short note we give an interpretation of the Kirchheim-Magnani 
counterexample to 
metric differentiability.  In fact we show that they construct 
something which fails shortly from being a  dilatation structure. 

Dilatation structures have been introduced in \cite{buligadil1}. These
structures are related to conical group 
\cite{buligacont}, which  form a particular class of contractible groups 
and are a slight generalization of Carnot groups. 

Carnot groups, in particular the Heisenberg group, appear as infinitesimal
models of sub-riemannian manifolds \cite{bell}, \cite{gromo}. In \cite{buligasr}
we explain how the formalism of dilatation structures applies to sub-riemannian
geometry. 

Further on we shall use the notations, definitions and results concerning 
dilatation structures, as found in \cite{buligadil1}, \cite{buligacont} or 
\cite{buligasr}. 

We shall construct a structure $(H(1), \rho,\bar{\delta})$ on  $H(1)$  which 
satisfies all axioms of a dilatation structure, excepting A3 and A4.  
 We prove  that for  $(H(1), \rho,\bar{\delta})$ the axiom  A4 implies A3. 
 Finally  we  prove that A4 for $(H(1), \rho,\bar{\delta})$ is equivalent with 
 $id$ metrically differentiable from $(H(1),d)$ to $(H(1),\rho)$, where 
 $d$ is a left invariant CC distance.  

For other relations between dilatation structures and differentiability in 
metric spaces see \cite{buligarn}.

\section{Metric differentiability for conical groups}
\label{kirma}

The general definition of metric differentiability for conical 
groups is formulated exactly as the same notion for Carnot groups. 

\begin{defi}
Let $(N,d,\delta)$ be a conical group. A continuous function $\eta:N\rightarrow [0,+\infty)$ is a 
seminorm if:
\begin{enumerate}
\item[(a)] $\displaystyle \eta(\delta_{\varepsilon} x) = \varepsilon \eta(x)$ for any $x\in N$ and $\varepsilon> 0$,
\item[(b)] $\eta(x y) \leq \eta(x)+\eta(y)$ for any $x,y\in N$.
\end{enumerate}
 Let $(N,\delta,d)$ be a conical group,  $(X,\rho)$ a metric space, $A\subset N$ an open set and 
  $x\in A$.    A function  $f:A \rightarrow X$ is metrically differentiable in $x$ if there is a seminorm 
 $\displaystyle \eta_{x}:N\rightarrow [0,+\infty)$ such that 
 $$\mid \frac{1}{\varepsilon} \rho(f(x\delta_{\varepsilon}v), f(x))  - \eta_{x}(v) \mid \rightarrow 0$$
 as $\varepsilon\rightarrow 0$, uniformly with respect to $v$ in compact neighbourhood of the neutral element $e\in N$. 
 \label{defmdif}
 \end{defi}

 \section{Kirchheim-Magnani counterexample to metric differentiability}

 For the elements of the Heisenberg group $\displaystyle H(1)= \mathbb{R}^{2}\times\mathbb{R}$  we use the notation $\displaystyle \tilde{x}=(x,\bar{x})$, with $\displaystyle \tilde{x} \in H(1), x\in \mathbb{R}^{2}, \bar{x}\in\mathbb{R}$. 
 In this subsection we shall use the following operation on $H(1)$: 
 $$\tilde{x}\tilde{y} = (x,\bar{x})(y,\bar{y}) = (x+y, \bar{x}+\bar{y}+ 2\omega(x,y)) , $$
 where $\omega$ is the canonical symplectic form on $\displaystyle \mathbb{R}^{2}$. 
 On $H(1)$ we consider the left invariant distance $d$ uniquely determined by the formula: 
 $$d((0,0), (x,\bar{x})) = \max\left\{ \|x\|, \sqrt{\mid\bar{x}\mid}\right\}  . $$
 
 The  construction by Kirchheim and Magnani is described further. Take an invertible, non decreasing function $g:[0,+\infty)\rightarrow [0,+\infty)$, continuous at $0$, such that $g(0)=0$. 
 
For a function $g$ which is well chosen,   the function $\rho:H(1)\rightarrow [0,+\infty)$, 
$$\rho(\tilde{x}) = \max\left\{ \|x\|, g(\mid\bar{x}\mid)\right\}$$
induces a left invariant invariant distance on $H(1)$ (we use the same symbol) 
$$\rho(\tilde{x}, \tilde{y}) = \rho(\tilde{x}^{-1} \tilde{y}) . $$
In order to check this it is sufficient to prove that for any $\displaystyle \tilde{x}, \tilde{y}\in H(1)$ we have 
$$\rho(\tilde{x}\tilde{y}) \leq \rho(\tilde{x}) + \rho(\tilde{y})  , $$
and that $\displaystyle  \rho(\tilde{x})= 0$ if and only if $\displaystyle  \tilde{x}=(0,0)$. 
The following result is  theorem 2.1 \cite{kirch}. 
 
 \begin{thm} (Kirchheim-Magnani) 
 If  the function $g$ has the expression  $$g^{-1}(t) = k(t) + t^{2}  $$
for any $t>0$, where $k:[0,+\infty)\rightarrow[0,+\infty)$ is a convex function, strictly increasing, continuous at $0$, and such that $k(0)=0$, then the function $\rho$ induces a left invariant distance (denoted also by 
$\rho$). Moreover, the identity function $id$ is 1-Lipschitz from   $(H(1),d)$ to $(H(1),\rho)$. 
\label{p410}
\end{thm}

 \section{Interpretation in terms of dilatation structures}
 
 Further we shall work with a function $g$ satisfying the hypothesis of theorem \ref{p410}, and with the
 associated function $\rho$ described in the previous subsection.

\begin{defi}
Define  for any $\varepsilon> 0$, the function 
 $$\bar{\delta}_{\varepsilon}(x,\bar{x}) = (\varepsilon x , sgn(\bar{x}) g^{-1}\left( \varepsilon g(\mid \bar{x}\mid)\right) ) $$
 for any $\displaystyle \tilde{x}=(x,\bar{x}) \in H(1)$. 
 
We define the following field of dilatations $\bar{\delta}$ by: for any $\varepsilon>0$ and $\tilde{x},\tilde{y}\in H(1)$ 
let 
$$\bar{\delta}^{\tilde{x}} \tilde{y} = \tilde{x} \bar{\delta}\left(\tilde{x}^{-1}\tilde{y}\right) \quad .$$
For any $\varepsilon>0$ and $\tilde{x},\tilde{y}\in H(1)$ we define 
$$\bar{\beta}_{\varepsilon}(\tilde{x},\tilde{y}) = \bar{\delta}_{\varepsilon^{-1}}\left( \bar{\delta}_{\varepsilon}(\tilde{x}) \bar{\delta}_{\varepsilon}(\tilde{y})\right) . $$

\end{defi}

We want to know when $(H(1),\rho,\bar{\delta})$ is a dilatation structure.

\begin{prop}
The  structure  $(H(1), \rho,\bar{\delta})$ satisfies the axioms A0, A1, A2. Moreover, A4 implies A3.
\label{p411}
\end{prop}

\paragraph{Proof.}
It is  easy to check that for any $\varepsilon, \mu \in (0,+\infty)$ we have 
$$\bar{\delta}_{\varepsilon}\bar{\delta}_{\mu} = \bar{\delta}_{\varepsilon\mu}  $$
and that $id = \delta_{1}$. 

Moreover, from $g$ non decreasing and continuous at $0$ we deduce that 
$$\lim_{\varepsilon\rightarrow 0} \bar{\delta}_{\varepsilon} \tilde{x} = (0,0)  , $$
uniformly with respect to $\displaystyle \tilde{x}$ in compact sets.

Another  computation shows that 
$$\rho(\bar{\delta}_{\varepsilon} \tilde{x}) = \varepsilon \rho(\tilde{x}) $$
for any $\displaystyle \tilde{x} \in H(1)$ and $\varepsilon> 0$.  
Otherwise stated, the function $\rho$ is homogeneous with respect to $\bar{\delta}$. 

All that is left to prove is that  A4 implies A3.
Remark that  $\bar{\delta}$ is left invariant (in the sense of transport by left translations in $H(1)$) and the distance $\rho$ is also left invariant. Then  axiom A4 takes the form: 
there exists the limit 
\begin{equation}
\lim_{\varepsilon\rightarrow 0} \bar{\beta}_{\varepsilon}(\tilde{x},\tilde{y}) = \bar{\beta}(\tilde{x},\tilde{y}) \in H(1)
\label{needpcnt}
\end{equation}
uniform with respect to $\tilde{x},\tilde{y}\in K$, $K$ compact set. 

From the homogeneity of the  function $\rho$  with respect to $\bar{\delta}$ we deduce that for any $\tilde{x},\tilde{y}\in H(1)$ we have: 
$$\frac{1}{\varepsilon} \rho\left( \bar{\delta}_{\varepsilon}(\tilde{x}), \bar{\delta}_{\varepsilon}(\tilde{y})\right) = \rho(\bar{\beta}_{\varepsilon}(\tilde{x}^{-1},\tilde{y})) . $$
From the left invariance of $\bar{\delta}$ and $\rho$ it follows that A4 implies A3. \quad $\square$

 \begin{thm}
 If the triple $(H(1),\rho,\bar{\delta})$ is a dilatation structure then $id$ is metrically differentiable from 
 $(H(1),d)$ to $(H(1),\rho)$. 
 \label{pcnt}
 \end{thm}
 
 \paragraph{Proof.} 
We know that  the triple $(H(1),\rho,\bar{\delta})$ is a dilatation structure if and only if (\ref{needpcnt}) 
is true. Taking (\ref{needpcnt}) as hypothesis we deduce that the identity function is derivable from 
$(H(1),d,\delta)$  to $(H(1),\rho,\bar{\delta})$. Indeed, computation shows that  $id$ derivable is 
equivalent to the existence of the limit 
$$\lim_{\varepsilon\rightarrow 0} \bar{\delta}_{\varepsilon^{-1}} \delta_{\varepsilon} \tilde{u} = 
(u, sgn(\bar{u}) g^{-1}\left( \lim_{\varepsilon\rightarrow 0} \frac{1}{\varepsilon} g(\varepsilon^{2} \mid\bar{u}\mid)\right)) $$
uniform with respect to $\tilde{u}$ in compact set. Therefore the function $id$ is derivable everywhere if 
and only if the uniform limit, with respect to $\bar{u}$ in compact set: 
\begin{equation}
A(\bar{u}) = \lim_{\varepsilon\rightarrow 0} \frac{1}{\varepsilon} g(\varepsilon^{2} \mid\bar{u}\mid) 
\label{ider}
\end{equation}
exists. We want to show that (\ref{needpcnt}) implies the existence of this limit. 

For this we shall use an equivalent (isomorphic) description of $(H(1),\rho,\bar{\delta})$. Consider the 
function $F:H(1)\rightarrow H(1)$, defined by 
$$F(x,\bar{x}) = (x, sgn(\bar{x}) g(\mid \bar{x}\mid) ) . $$
The function $F$ is invertible because $g$ is invertible. For any $\varepsilon>0$ let 
$\displaystyle \hat{\delta}_{\varepsilon}$ be the usual dilatations: 
$$\hat{\delta}_{\varepsilon} (x,\bar{x}) = (\varepsilon x, \varepsilon \bar{x}) .$$
It is then straightforward that 
$$\bar{\delta}_{\varepsilon} = F^{-1} \hat{\delta}_{\varepsilon} F , $$
for any $\varepsilon>0$.

The function $F$ can be made into a group isomorphism by re-defining the group operation on $H(1)$ 
$$\tilde{x}\cdot\tilde{y} = F( (x,h(\bar{x}))(y,h(\bar{y}))  , $$
where $h$ is the function 
$$h(t) = sgn(t) ( t^{2} + k(\mid t \mid)) . $$
Let $\mu$ be the transported left invariant distance on $H(1)$, defined by 
$$\mu(F(\tilde{x}), F(\tilde{y})) = \rho(\tilde{x},\tilde{y}) . $$
Remark that $\mu$ has the simple expression 
$$ \mu((0,0),(x,\bar{x})) =   \max \left\{ \mid x\mid, \mid \bar{x}\mid \right\} . $$
Exactly as before we can construct the structure 
$\hat{\delta}$ by 
$$\hat{\delta}^{\tilde{x}}_{\varepsilon} \tilde{y} = \tilde{x} \cdot \hat{\delta}_{\varepsilon} \left( \tilde{x}^{-1}\cdot \tilde{y}\right)    .  $$
We get a dilatation structure  $(H(1),\mu,\hat{\delta})$ isomorphic with $(H(1),\rho,\bar{\delta})$. 

The identity function $id$ is derivable from  $(H(1),d,\delta)$  to $(H(1),\rho,\bar{\delta})$ if and only 
if the function $F$ is derivable from $(H(1),d,\delta)$  to $(H(1),\mu,\hat{\delta})$. 

 The axiom A4 for the dilatation structure  $(H(1),\mu,\hat{\delta})$ implies that   for any 
 $\tilde{x},\tilde{y} \in H(1)$ the limit exists 
$$\lim_{\varepsilon\rightarrow 0} \frac{1}{\varepsilon} g \left( \mid \varepsilon^{2}\left( \frac{1}{2}\omega(x,y) + \mid \bar{x}\mid  \bar{x} + \mid \bar{x}\mid  \bar{x} \right) + sgn(\bar{x}) k(\varepsilon\mid \bar{x}\mid ) + sgn(\bar{y}) k(\varepsilon\mid \bar{y}\mid ) \mid \right) , $$
uniform with respect to $\tilde{y}$ in compact set.   Take in the previous limit $\bar{x}=\bar{y}=0$ and 
denote $\displaystyle \bar{u} = \frac{1}{2}\omega(x,y)$. We get (\ref{ider}), therefore we proved that 
$id$ is derivable  from  $(H(1),d,\delta)$  to $(H(1),\rho,\bar{\delta})$. 

Finally, the derivability of $id$ implies the metric differentiability. Indeed, we use (\ref{ider}) to 
compute   $\nu$, the metric differential of $id$. We obtain that 
$$\nu_{\tilde{x}} = \mu((x,A(\bar{x}))) = \max\left\{ \mid x \mid,  A(\bar{u}) \right\}  \quad . $$
The proof is done. \quad $\square$

In the counterexample of Kirchheim and Magnani the identity function $id$ is not metric differentiable, 
therefore the corresponding triple $(H(1),\rho,\bar{\delta})$ is not a dilatation structure.

\end{document}